\documentclass[letterpaper, 10 pt, conference]{ieeeconf}
\IEEEoverridecommandlockouts
\usepackage{graphics} 
\usepackage{epsfig} 
\usepackage{mathptmx} 
\usepackage{times} 
\usepackage{amsmath} 
\usepackage{amssymb}  

\usepackage{subfigure}
\usepackage{algorithm}
\usepackage{algpseudocode}

\newtheorem{theo}{Theorem}
\newtheorem{lemma}{Lemma}
\newtheorem{assum}{Assumption}

\title{\LARGE \bf
Infinite-horizon Risk-constrained Linear Quadratic Regulator  with Average Cost
}

\author{Feiran Zhao, Keyou You, Tamer~Ba\c{s}ar
\thanks{This research was supported by National Natural Science Foundation of China under Grant no. 62033006.}
\thanks{F. Zhao and K. You are with the Department of Automation and BNRist, Tsinghua University, Beijing 100084, China. e-mail: zhaofr18@mails.tsinghua.edu.cn, youky@tsinghua.edu.cn.}%
\thanks{T.~Ba\c{s}ar is with the Coordinated Science Laboratory, University of Illinois at Urbana-Champaign, Urbana, IL 61801 USA. e-mail:basar1@illinois.edu.}
}

\begin{document}

\maketitle
\thispagestyle{empty}
\pagestyle{empty}

\begin{abstract}
The behaviour of a stochastic dynamical system may be largely influenced by those low-probability, yet extreme events. To address such occurrences, this paper proposes an infinite-horizon risk-constrained Linear Quadratic Regulator (LQR) framework with time-average cost. In addition to the standard LQR objective, the average one-stage predictive variance of the state penalty is constrained to lie within a user-specified level. By leveraging the duality, its optimal solution  is first shown to be stationary and affine in the state, i.e., $u(x,\lambda^*) = -K(\lambda^*)x + l(\lambda^*)$, where $\lambda^*$ is an optimal multiplier, used to address the risk constraint. Then, we establish the stability of the resulting closed-loop system. Furthermore, we propose a primal-dual method with sublinear convergence rate to find an optimal policy $u(x,\lambda^*)$. Finally, a numerical example is provided to demonstrate the effectiveness of the proposed framework and the primal-dual method.
\end{abstract}

\section{Introduction}
Stochastic optimal control is a long-studied framework for the dynamical systems with uncertain variables~\cite{aastrom2012introduction,bertsekas1995dynamic}. For example, the Linear Quadratic Regulator (LQR) with noisy input involves minimization of time-average cumulative quadratic cost in the expectation, which is known to be risk-neutral~\cite{whittle1981risk, anderson2007optimal}. That is, its performance may significantly degrade due to the presence of the low-probability, yet extreme external noises. For those safety-critical applications e.g., the autonomous vehicles, the non-typical events can even lead to catastrophic consequences. Thus, a fundamental problem is to address the potential risk in stochastic systems.

For decades, risk-aware control has drawn an increasing research interest for its promise to deal with unexpected noises~\cite{huang2017risk, bauerle2014more, sopasakis2019riskECC, sopasakis2019riskA, jiang2018risk}. Based on the optimal control framework, it typically compensates for the risk by considering an exponentiation of the regulation cost~\cite{whittle1981risk, moore1997risk,ito2018risk,speyer1992optimal}. However, the noise distribution is typically limited to a specific class to render a well-defined optimization problem. As a consequence, some risk-aware controllers are unable to handle noises with asymmetric structures e.g., skewed distributions. A well-known instance is the Linear Exponential Quadratic Gaussian (LEQG) control, where the exponential cost is interpreted as a linear combination of quadratic cost and its variance and higher moments. While it yields a simple closed-form controller, the process noise is assumed to be Gaussian with zero mean. In contrast to the exponential approach, the risk awareness to heavy-tailed distributions can be achieved by optimizing a risk measure~\cite{borkar2014risk, chapman2019risk, shapiro2014lectures,  di2012policy,tessler2018reward,chow2017risk} e.g., Conditional Value at Risk (CVaR)~\cite{rockafellar2000optimization}. However, it is challenging to obtain a closed-form solution, and approximations are widely used for tractability. Recently, a new risk measure for LQR control has been introduced, i.e., the cumulative expected one-step predictive variance of the state penalty~\cite{tsiamis2020risk}. By setting it as a constraint,~\cite{tsiamis2020risk} has further proposed a finite-horizon risk-constrained LQR framework. Under mild statistical conditions on the noises, it leads to a non-stationary closed-form controller which is affine in the state.

This work can be viewed as an extension of the discrete-time risk-constrained LQR~\cite{tsiamis2020risk} to the infinite-horizon setting. Specifically, we aim to find a control sequence to minimize a time-average LQR cost subject to an average one-stage predictive variance constraint of the state penalty. We note that this extension is non-trivial in at least two aspects. First, the finite-horizon risk-constrained LQR problem can be converted to a Quadratically Constrained Quadratic Program (QCQP), which is essentially convex. In sharp contrast, there are infinite number of variables in our constrained optimization problem, and it cannot be handled via functional analysis due to the limit in the average cost formulation. Second, a solution obtained by simply letting the horizon tend to infinity, though plausible, may not be optimal. In fact, we rigorously prove its optimality by establishing an Average-Cost Optimality Equation (ACOE). Our contribution lies in addressing satisfactorily the above issues, and further showing by duality that an optimal solution is stationary and affine in the state, i.e., $u(x, \lambda^*) = -K(\lambda^*)x + l(\lambda^*)$, where $\lambda^*$ is the optimal multiplier used to address the risk constraint. Moreover, we propose a primal-dual method with a sublinear convergence rate to find an optimal policy $u(x, \lambda^*)$. As a comparison, \cite{tsiamis2020risk} applies simple bisection to search an optimal multiplier, yet without convergence analysis.

Our work is pertinent to LQR with cumulative cost constraints~\cite{lim1999stochastic, lim1999linear, lim1997quasi, lim1996linearly, bakolas2016optimal, tsiamis2020risk}. In physical-world applications, many design objectives can be expressed by quadratic cost constraints. In fact, the proposed risk constraint is also shown to be quadratic in the state. The finite-horizon constrained LQR problem has been solved in both discrete-time~\cite{bakolas2016optimal, tsiamis2020risk} and continuous-time~\cite{lim1999stochastic, lim1999linear, lim1997quasi, lim1996linearly} settings by leveraging the convexity, leading to a simple non-stationary feedback controller. To the best of our knowledge, however, there are no such results in the infinite-horizon setting since infinite-dimensional stochastic optimization with constraints is generally difficult.

The remainder of this paper is organized as follows. In Section \ref{sec:prob}, we introduce the infinite-horizon risk-constrained LQR problem with time-average cost. In Section \ref{sec:opt}, we first reformulate the risk constraint as a time-average cost that is quadratic in the state. Then, we establish an optimality condition for the infinite-horizon risk-constrained LQR by duality. In Section \ref{sec:pd}, we propose a primal-dual method with convergence guarantees to find an optimal policy. In Section \ref{sec:exp}, we validate our results via simulations.

\section{Problem Formulation}\label{sec:prob}

For the LQR problem with external noises, we consider a discrete-time linear stochastic system with full state observations
\begin{equation}\label{equ:sys}
x_{t+1} = Ax_t + Bu_t+w_t,
\end{equation}
where $x_t \in \mathbb{R}^n$ denotes the state, $u_t \in \mathbb{R}^m$ is the control input and $w_t \in \mathbb{R}^d$ is the uncorrelated random noise. $A\in \mathbb{R}^{n \times n}$ and $B \in \mathbb{R}^{n \times m}$ are the model parameters.

The infinite-horizon LQR targets to find a control sequence $u = \{u_0, u_1, \dots\}$ in the form of $u_t = \pi_t(h_t)$ with the system history trajectory $h_t = \{x_0,u_0,\cdots, x_{t-1},u_{t-1}, x_t\}$, to minimize the following time-average cost, i.e.,
\begin{equation}\label{prob:lqr}
\begin{aligned}
\text { minimize } &~~ \limsup\limits_{T \rightarrow  \infty}  \frac{1}{T} \mathbb{E}  \sum_{t=0}^{T-1}(x_{t}^{\top} Q x_{t}+u_{t}^{\top} R u_{t})\\
\text {subject to} &~~(\ref{equ:sys})
\end{aligned}
\end{equation}
where the expectation is taken with respect to the random noise $\{w_t\}$ and the policy $\{\pi_t\}$, which is not necessarily deterministic. Throughout the paper, we make the following assumption standard in control theory~\cite{bertsekas1995dynamic}.

\begin{assum}
	\label{assumption}
	$Q$ is positive semi-definite and $R$ is positive definite. The pair $(A,B)$ is stabilizable and $(A,Q^{{1}/{2}})$ is observable.
\end{assum}

Under Assumption \ref{assumption}, a unique optimal policy to (\ref{prob:lqr}) is known to be stationary and linear in the state when $w_t$ has zero mean, i.e., $u_t=-Kx_t$. Under such a policy, however, the system state may be significantly affected by the extreme events as the LQR only minimizes the expected cost.

In this work, we extend the finite-horizon risk-constrained LQR in~\cite{tsiamis2020risk} to the infinite-horizon setting. Specifically, we aim to minimize the average cost subject to a one-step predictive state variability constraint, i.e.,
\begin{equation}\label{equ:rclqr}
\begin{aligned}
\text {minimize } &~ \limsup\limits_{T \rightarrow  \infty} \frac{1}{T} \mathbb{E}\sum_{t=0}^{T-1}(x_{t}^{\top} Q x_{t}+u_{t}^{\top} R u_{t})\\
\text {subject to}& ~~(\ref{equ:sys}) ~\text{and} \\
&\limsup\limits_{T \rightarrow  \infty} \frac{1}{T}  \mathbb{E} \sum_{t=0}^{T-1}(x_{t}^{\top} Q x_{t}-\mathbb{E}[x_{t}^{\top} Q x_{t}|h_t])^{2} \leq \rho \\
\end{aligned}
\end{equation}
where $\rho > 0$ is a user-defined constant for the risk tolerance. In contrast to the LEQG~\cite{whittle1981risk}, we do not require the noise $w_t$ to be Gaussian. Instead, we only assume that $w_t$ has a finite fourth-order moment~\cite{tsiamis2020risk}. Note that the limit of the average cost may not exist as $T \rightarrow \infty$, thus a supremum must be taken in both the objective and the constraint.

In this paper, we show via duality that an optimal controller to (\ref{equ:rclqr}) is stationary and affine in the state, i.e., $u(x, \lambda^*) = - K(\lambda^*)x +l(\lambda^*)$, where $\lambda^*$ is an optimal multiplier to address the risk constraint. Then, we prove the stability of the resulting closed-loop system. Furthermore, we propose a primal-dual method with sublinear convergence rate to find an optimal policy $u(x, \lambda^*)$.

\section{An Optimal Controller to the Infinite-horizon Risk-constrained LQR}\label{sec:opt}
In this section, we first reformulate (\ref{equ:rclqr}) as a quadratically  constrained quadratic program (QCQP) problem. Then, by exploiting properties of the Lagrangian function, we find an optimal controller that solves (\ref{equ:rclqr}), which is also able to stabilize the system.

\subsection{Reformulation of (\ref{equ:rclqr})}
Define the mean $\bar{w} = \mathbb{E}[w_t]$, the covariance $W = \mathbb{E}[(w_t - \bar{w})(w_t - \bar{w})^{\top}]$ and other higher-order weighted statistics of $w_t$ by
\begin{align*}
M_{3} &= \mathbb{E}[(w_{t}-\bar{w})(w_{t}-\bar{w})^{\top} Q(w_{t}-\bar{w})], \\
m_{4} &= \mathbb{E}[(w_{t}-\bar{w})^{\top} Q(w_{t}-\bar{w})-\operatorname{tr}\{WQ\}]^{2}. \label{fouth}
\end{align*}
Then, by~\cite{tsiamis2020risk}, we can reformulate (\ref{equ:rclqr}) as
\begin{equation}\label{prob:new_rclqr}
\begin{aligned}
\mathop{\text {minimize }}\limits_{u_0,u_1,\dots} &~ J(u):= \limsup\limits_{T \rightarrow  \infty} \frac{1}{T} \mathbb{E}\sum_{t=0}^{T-1}(x_{t}^{\top} Q x_{t}+u_{t}^{\top} R u_{t})\\
\text {subject to}& ~~(\ref{equ:sys}) ~\text{and} \\
&J_c(u) = \limsup\limits_{T \rightarrow  \infty} \frac{1}{T}  \mathbb{E} \sum_{t=0}^{T-1}(4x_{t}^{\top} QWQ x_{t} + 4x_{t}^{\top}QM_3)\\
& \leq \bar{\rho}= \rho - m_4 + 4\operatorname{tr}\{(WQ)^2\}.
\end{aligned}
\end{equation}
where the constraint is quadratic in the state. In contrast to its finite-horizon setting~\cite{tsiamis2020risk}, we have here an infinite number of optimization variables $\{u_0, u_1, \dots\}$ in (\ref{prob:new_rclqr}). It is not possible to approach the average cost constrained problem (\ref{prob:new_rclqr}) via the functional analysis. Moreover, we cannot directly apply the dynamic programming paradigm to minimize $J(u)$ due to the presence of the risk constraint. In the rest of this section, we approach (\ref{prob:new_rclqr}) by building intuitions from duality theory.

Let $\lambda \geq 0$ denote the multiplier associated with (\ref{prob:new_rclqr}), $Q_{\lambda} = Q+4 \lambda Q W Q$ and $S = 2 QM_3$. Define the Lagrangian as
\begin{equation}\label{def:L}
\begin{aligned}
L(u,\lambda) &= J(u)+\lambda (J_c(u) - \bar{\rho})\\
&= \limsup\limits_{T \rightarrow  \infty} \frac{1}{T}  \mathbb{E} \sum_{t=0}^{T-1}(x_{t}^{\top} Q_{\lambda} x_{t}+ 2\lambda S^{\top}x_t +u_{t}^{\top} R u_{t}) -\lambda \rho
\end{aligned}
\end{equation}
and consider the following optimization problem
\begin{equation}\label{prob:lag}
\text{minimize}_{u} ~L(u, \lambda).
\end{equation}
We show that (a) for a given $\lambda$, solving (\ref{prob:lag}) yields a stationary optimal policy $u(x, \lambda)$ and (b) there exists an optimal $\lambda^*$ such that $u(x, \lambda^*)$ is a solution to the original problem (\ref{prob:new_rclqr}).

\subsection{Optimality Equation for (\ref{prob:lag})}

%

We first show that an optimal solution to (\ref{prob:lag}) is stationary and affine, i.e., $u(x, \lambda) = - K(\lambda) x + l(\lambda)$. To this end, we build insights by considering the finite-horizon cost
\begin{equation}\label{def:finite_cost}
\mathbb{E} \sum_{t=0}^{T-1}(x_{t}^{\top} Q_{\lambda} x_{t}+ 2\lambda S^{\top}x_t +u_{t}^{\top} R u_{t}) .
\end{equation}
We reorganize the main results in \cite{tsiamis2020risk} in the following lemma. Note that we use $u$ and $\pi$ interchangeably to denote a policy.
\begin{lemma}(\cite[Theorem 2]{tsiamis2020risk})\label{lem:value}
	For a given $\lambda$, an optimal value of (\ref{def:finite_cost}) has a quadratic form
	$$V_t(x_t) = x_t^{\top}P_tx_t + g_tx_t + z_t, ~~\forall t \in \{0,1,\dots,T\},$$
	and an optimal policy to (\ref{prob:lag}) is
	\begin{equation}
	u_t = -K_tx_t + l_t, ~~\forall t \in \{0,1,\dots,T-1\},
	\end{equation}
	where for $t \in \{1,\dots,T\}$,
	\begin{align*}
	P_{t-1} &= Q_{\lambda} + A^{\top}P_{t}A - A^{\top}P_{t}B(R+B^{\top}P_{t}B)^{-1}B^{\top}P_{t}A\\
	K_{t-1} &= - (R + B^{\top}P_{t-1}B)^{-1}B^{\top}P_{t-1}Ax  \\
	g_{t-1}^{\top} &= (2\bar{w}^{\top}P_{t} + g_{t}^{\top})(A-BK_{t-1}) + 2\lambda S^{\top},\\
	l_{t-1} &= - \frac{1}{2}(R + B^{\top}P_{t-1}B)^{-1}B^{\top}(2P_{t-1}\bar{w} + g_{t-1})\\
	z_{t-1} &= z_{t} + \text{tr}\{P_{t}(W+\bar{w}\bar{w}^{\top})\} +g_{t}^{\top}\bar{w}  +l_{t-1}(R + B^{\top}P_{t}B)^{-1} l_{t-1}.
	\end{align*}
	with terminal values $P_T = 0$, $g_T = 0$ and $z_T = 0$.
\end{lemma}

By Lemma \ref{lem:value}, the optimal value of the average $T$-stage cost
\begin{equation}\label{equ:average}
\frac{1}{T}  \mathbb{E} \sum_{t=0}^{T-1}(x_{t}^{\top} Q_{\lambda} x_{t}+ 2\lambda S^{\top}x_t +u_{t}^{\top} R u_{t}) -\lambda \rho
\end{equation}
is equal to
\begin{align*}
&\frac{1}{T}(x_0^{\top}P_0x_0 + g_0^{\top} + \sum_{t=0}^{T-1}(\text{tr}\{P_{t+1}(W+\bar{w}\bar{w}^{\top})\} +g_{t+1}^{\top}\bar{w} \\
&~~~ +l_t(R + B^{\top}P_{t+1}B)^{-1} l_t)) - \lambda \rho.
\end{align*}

Clearly, letting $T \rightarrow \infty$, $P_t$ yields in the limit a positive semi-definite matrix $P$, which is given by the solution of the algebraic Riccati equation
\begin{equation}\label{equ:riccati}
P = Q_{\lambda} + A^{\top}PA - A^{\top}PB(R+B^{\top}PB)^{-1}B^{\top}PA.
\end{equation}

Hence, the control gain $K_t$ converges to
\begin{equation}\label{k}
K(\lambda) = - (R + B^{\top}PB)^{-1}B^{\top}PA
\end{equation}
which is also able to stabilize the closed-loop system~\cite{bertsekas1995dynamic}, i.e., $\rho(A-BK)<1$. Similarly, it follows that $g_t$ converges to a fixed point given by
$$g^{\top} = (2\bar{w}^{\top}P + g^{\top})(A-BK) + 2\lambda S^{\top}$$
and $l_t$ converges to
\begin{equation}\label{l}
l(\lambda) = - \frac{1}{2}(R + B^{\top}PB)^{-1}B^{\top}(2P\bar{w} + g).
\end{equation}
Furthermore, as $T \rightarrow \infty$, the average cost (\ref{equ:average}) tends to
\begin{equation}\label{h}
h(\lambda) = \text{tr}\{P(W+\bar{w}\bar{w}^{\top})\} +g^{\top}\bar{w} +l(R + B^{\top}PB)^{-1} l - \lambda \rho.
\end{equation}

We rigorously prove in the following theorem that the stationary policy $u(x, \lambda) = -K(\lambda)x + l(\lambda)$ is indeed an optimal solution to (\ref{prob:lag}) by establishing an average-cost optimality equation (ACOE).

\begin{theo}\label{theo:average}
	For a fixed $\lambda$, an optimal policy to (\ref{prob:lag}) is stationary and affine in the state, i.e.,
	\begin{equation}\label{equ:opt_policy}
	u(x, \lambda) = -K(\lambda)x + l(\lambda).
	\end{equation}
	where $K(\lambda)$ and $l(\lambda)$ are given by (\ref{k}) and (\ref{l}), respectively. Moreover, $h(\lambda)$ in (\ref{h}) is the optimal value of the Lagrangian $L(u,\lambda)$.
\end{theo}

\begin{proof}
	By the definition of $h(\lambda)$, $P$, $g$ and $u(x, \lambda)$, it is straightforward to show that the following ACOE holds
	\begin{align*}
	&h(\lambda) + x^{\top}Px + g^{\top}x \\
	&= \min_{u} \mathbb{E}[x^{\top} Q_{\lambda} x+ 2\lambda S^{\top}x +u^{\top} Ru +g^{\top}(Ax + Bu+ w)\\
	&~~~+ (Ax + Bu+w)^{\top}P(Ax + Bu + w)],
	\end{align*}
	where the minimum of the right hand side (RHS) is attained at $u(x, \lambda)$.
	
	We show that the $h(\lambda)$ is the optimal value of $L(u,\lambda)$. It follows from the ACOE that for any policy $\pi$,
	\begin{align*}
	&x_0^{\top}Px_0 + g^{\top}x_0 \\
	&\leq \mathbb{E}[x_0^{\top} Q_{\lambda} x_0+ 2\lambda S^{\top}x_0 +u_0^{\top} Ru_0 + x_1^{\top}Px_1 + g^{\top} x_1 | x_0, \pi] \\
    &~~~- h(\lambda) \\
	&\leq \dots \\
	& \leq \mathbb{E}[\sum_{t=0}^{T-1}(x_t^{\top} Q_{\lambda} x_t+ 2\lambda S^{\top}x_t +u_t^{\top} Ru_t) | x_0, \pi] \\
	&~~~+ \mathbb{E}[x_T^{\top}Px_T + g^{\top}x_T | x_0, \pi] - Th(\lambda)
	\end{align*}
	
	Dividing by $T$, we have
	\begin{equation}\label{equ:opt}
	\begin{aligned}
	h(\lambda) &\leq -\frac{1}{T}(x_0^{\top}Px_0 + g^{\top}x_0) + \frac{1}{T}\mathbb{E}[x_T^{\top}Px_T + g^{\top}x_T]\\
	&~~~+ \frac{1}{T} \mathbb{E}[\sum_{t=0}^{T-1}(x_t^{\top} Q_{\lambda} x_t+ 2\lambda S^{\top}x_t +u_t^{\top} Ru_t) | x_0, \pi].
	\end{aligned}
	\end{equation}
	
	Since for a policy $\pi$ that satisfies $\sup_{t \geq 0}\mathbb{E}\|x_t\|^2 < \infty$, it follows that
	$$\lim_{T \rightarrow \infty} \frac{1}{T}\mathbb{E}\|x_T^{\top}Px_T + g^{\top}x_T\| = 0, $$
	by limiting $T \rightarrow \infty$, it follows from (\ref{equ:opt})  that
	$$
	h(\lambda) \leq \lim_{T \rightarrow \infty} \sup \frac{1}{T} \mathbb{E}[\sum_{t=0}^{T-1}(x_t^{\top} Q_{\lambda} x_t+ 2\lambda S^{\top}x_t +u_t^{\top} Ru_t) | x_0, \pi],
	$$
	which completes the proof.
\end{proof}

The stability of the resulting closed-loop system follows from standard control theory~\cite{anderson2007optimal}.
\begin{lemma}
	For a given $\lambda$, the stationary policy $u(x, \lambda)$ in (\ref{equ:opt_policy}) is able to stabilize the system, i.e., $\rho(A - BK(\lambda))<1$.
\end{lemma}
\begin{proof}
	Since $K(\lambda) = - (R + B^{\top}PB)^{-1}B^{\top}PA$ where $P$ is given by the Riccati equation (\ref{equ:riccati}), it follows from~\cite{anderson2007optimal} that $\rho(A - BK(\lambda))<1$.
\end{proof}

\subsection{Optimality Conditions for (\ref{prob:new_rclqr})}
Define
\begin{equation}\label{def:mu}
\lambda^{*} = \inf \{\lambda \geq 0| J_c(u(x, \lambda)) \leq \rho\}.
\end{equation}
In this subsection, we show that $\lambda^{*}$ in (\ref{def:mu}) exists and $u(x, \lambda^*)$ is an optimal solution to (\ref{prob:new_rclqr}). 

\begin{lemma}\label{lem:continu}
	The risk constraint $J_c(u(x, \lambda))$ under policy $u(x, \lambda)$ is continuous in $\lambda$. 
\end{lemma}
\begin{proof}
	Clearly, the optimal policy $u(x, \lambda)$ is continuous with respect to $\lambda \geq 0$ due to the invertibility of $R + B^{\top}PB$. Hence, the risk constraint $J_c(u(x, \lambda))$ is continuous in $\lambda$.
\end{proof}

\begin{theo}\label{theo:duality}
	Suppose that Slater's condition holds, i.e., there exists a policy $\tilde{u}$ such that $J_c(\tilde{u}) < \rho$. Then,
	$$u(x, \lambda^*) = -K(\lambda^*)x + l(\lambda^*)$$
	 is an optimal solution to (\ref{prob:new_rclqr}).
\end{theo}

\begin{proof}
	We first show that (a) $\lambda^*$ defined in (\ref{def:mu}) exists and (b) the policy $u(x, \lambda^*)$ satisfies
	\begin{equation}\label{equ:saddle}
	\lambda^*(J_c(u(x, \lambda^*))- \rho)=0.
	\end{equation}
	
	(a) By the Slater's condition, there exists a constant $a>0$ such that $J_c(\tilde{u})+a \leq \rho$. We prove (a) by contradiction.
	
	Suppose that for all $\lambda \geq 0$, we have $J_c(u(x, \lambda)) > \rho$. Then
	\begin{align*}
	J(\tilde{u}) &\geq \min_{u}L(u,\lambda) - \lambda (J_c(\tilde{u}) - \rho)\\
	&\geq J(u(x, \lambda)) + \lambda J_c(u(x, \lambda)) -\lambda \rho - \lambda (J_c(\tilde{u}) - \rho)\\
	&\geq J(u(x, \lambda)) + \lambda (J_c(u(x, \lambda)) - J_c(\tilde{u})) \\
	&\geq J(u(x, \lambda)) + \lambda (J_c(u(x, \lambda)) - \rho + a) \\
	&> J(u(x, \lambda)) + \lambda a
	\end{align*}
	Let $\lambda \rightarrow \infty$, then $J(\tilde{u}) > \infty$, which contradicts the Slater's condition. Thus, $\lambda^*$ defined in (\ref{def:mu}) exists.
	
	(b) To show that $\lambda^*(J_c(u(x, \lambda^*))- \rho)=0$, we consider two cases. If $\lambda^* = 0$, then it trivially holds; otherwise, we must have $J_c(u(x, \lambda^*))=\rho$. Assuming that $\lambda^* >0$, it follows that $J_c(u(x, 0)) > \rho$. Since $\lambda^*$ is finite, there exists a multiplier $\lambda'$ such that $J_c(u(x, \lambda')) \leq \rho$. The continuity in Lemma leads to that $J_c(u(x, \lambda^*))=\rho$. Thus, $u(x, \lambda^*)$ with $\lambda^*$ given in (\ref{def:mu}) satisfies (\ref{equ:saddle}).
	
	Then, we prove the optimality of $u(x, \lambda^*)$. Denote $J^*$ as the optimal values of (\ref{equ:rclqr}). We have the following relations,
	\begin{equation}\label{equ:strong}
	\begin{aligned}
	J^* &\leq J(u(x, \lambda^*)) \\
	&= J(u(x, \lambda^*)) + \lambda^*(J_c(u(x, \lambda^*))- \rho) \\
	&= \min_{u} L(u, \lambda^*) \\
	&\leq \max_{\lambda} \min_{u} L(u, \lambda)
	\end{aligned}
	\end{equation}
	where the first equality follows from (\ref{equ:saddle}).
	
	Using the weak duality theorem~\cite{bertsekas1997nonlinear}, we have $J^* \geq \max_{\lambda} \min_{u} L(u, \lambda)$. Thus, the equality holds  throughout (\ref{equ:strong}), which implies that $u(x, \lambda^*)$ is an optimal policy.
\end{proof}

\section{Primal-dual Method to Solve the Risk-constrained LQR}\label{sec:pd}
In this section, we propose a primal-dual method with sublinear convergence rate to solve (\ref{prob:new_rclqr}).

By Theorem \ref{theo:duality}, there is no duality gap for (\ref{prob:new_rclqr}).  Thus, we can alternatively to solve the following dual problem of (\ref{prob:new_rclqr})
\begin{equation}\label{prob:dual}
\max_{\lambda\geq 0} D(\lambda) = \max_{\lambda\geq 0} \min_{u} L(u, \lambda),
\end{equation}
which is always concave in $\lambda$. By the dual theory~\cite{nesterov2013introductory, nedic2009subgradient, boyd2004convex}, a subgradient of $D(\lambda)$ is given as
\begin{equation}\label{def:sub}
d^k = J_c(u(x, \lambda^k)) - \bar{\rho},
\end{equation}
where $J_c(u(x, \lambda^k))$ is explicitly computed by the following lemma.
\begin{lemma}\label{lem:risk}
	For a stabilizing policy $u = -Kx + l$, we have
	$$
	J_c(u) = \mathrm{tr}\{P_c(W + (Bl+\bar{w})(Bl+\bar{w})^{\top})\} + g_c^{\top}(Bl+\bar{w}).
	$$
	where $P_c > 0$ is a unique solution of the Lyapunov equation
	$$
	P_{c} = 4QWQ + (A-B K)^{\top} P_{c}(A-B K),
	$$
	and
	$ g_{c}^{\top} = 2((Bl+\bar{w})^{\top}P_c(A-BK)+ 2M_3^{\top}Q)(I-A+BK)^{-1}.$
\end{lemma}

\begin{proof}
	Clearly, the risk measure $J_c(u)$ is finite under a stabilizing policy. Define the relative value function of the risk constraint by
	\begin{align*}\label{def:value_risk}
	V_c(x)&=\mathbb{E} \sum_{t = 0}^{\infty} [4x_{t}^{\top} QWQ x_{t} + 4x_{t}^{\top}QM_3 \\
	&~~~~~~~~~~- J_c(u)| x_{0}=x, u_t = -Kx_t +l].
	\end{align*}
	
	By using backward dynamic programming \cite{bertsekas1995dynamic}, it can be easily shown that $V_c(x)$ has a quadratic form, i.e.,
	$
	V_c(x)= x^{\top}P_cx + g_c^{\top}x + z_c,
	$
	where $P_c,  g_c, z_c$ are  to be determined.
	
	By Bellman equation~\cite{bertsekas1995dynamic}, it holds that
	\begin{align*}
	&x_t^{\top}P_cx_t + g_c^{\top}x_t + z_c\\
	& = 4x_t^{\top}QWQx_t + 4M_3^{\top}Qx_t  \\
	&~~~+ \mathbb{E}  [(A-BK)x_t +Bl + w_t]^{\top}P_c [(A-BK)x_t +Bl + w_t] \\
	&~~~ - J_c(u) + \mathbb{E}[g_c^{\top}((A-BK)x_t +Bl + w_t)] + z_c\\
	&= x_t^{\top}[4QWQ + (A-B K)^{\top} P_{c}(A-B K)]x_t \\
	&~~~+ [2(Bl+\bar{w})^{\top}P_c(A-BK) +4M_3^{\top}Q+ g_c^{\top}(A-BK)]x_t \\
	&~~~+ \mathrm{tr}\bigl[P_c(W + (Bl+\bar{w})(Bl+\bar{w})^{\top} )\bigr] + g_c^{\top}(Bl+\bar{w}) \\
	&~~~- J_c(u)  + z_c
	\end{align*}
	
	Noting that the equality holds for all $x_t \in \mathbb{R}^n$, we can only have
	$$
	J_c(u) = \mathrm{tr}\{P_c(W + (Bl+\bar{w})(Bl+\bar{w})^{\top})\} + g_c^{\top}(Bl+\bar{w}).
	$$
\end{proof}

We present our primal-dual method for (\ref{prob:dual}) in Algorithm \ref{alg}. Since $u(x, \lambda^k)$ is always able to stabilize the system, the subgradient $d^k$ and  $\lambda^{k}$ are bounded by some positive constants, i.e., $\|d^k\| \leq b$ and $\|\lambda^{k}\| \leq e$. The global convergence guarantee for Algorithm \ref{alg} then follows from the concavity of $D(\lambda)$. Denote $\max_{\lambda \geq 0}~ D(\lambda)$ in (\ref{prob:dual}) as $D^*$.

\begin{algorithm}[t]
	\caption{The primal-dual algorithm for the risk-constrained LQR}
	\label{alg}
	\begin{algorithmic}[1]
		\Require A multiplier $\lambda^1\geq 0$, and a set of stepsizes $\{\zeta^k\}$.
		\For{$k=1,2,\dots$}
		\State Solve $u(x, \lambda^k) = \text{argmin}_{u}~ L(u,\lambda^k)$ by Theorem \ref{theo:average}.
		\State Compute a subgradient $d^k$ by (\ref{def:sub}) and Lemma \ref{lem:risk}.
		\State Update the multiplier by $\lambda^{k+1} = [\lambda^{k} + \zeta^k \cdot d^k]_{+}$.
		\EndFor	
	\end{algorithmic}
\end{algorithm}

\begin{theo}\label{theo:model-based}
	Let $\bar{\lambda}^k = \frac{1}{k} \sum_{i=1}^{k} \lambda^{i}$. For $\zeta^k = \frac{1}{be}\sqrt \frac{2}{k}, k \in \{1,2,\dots\}$,  Algorithm \ref{alg} satisfies
	$$
	D^*-	D(\bar{\lambda}^k) \leq \frac{3be}{\sqrt{k}}.
	$$
\end{theo}
\begin{proof}
	By the definition of projection, it follows that
	\begin{align*}
	\|\lambda^{i+1} - \lambda^*\|^2 &\leq \| \lambda^{i} - \lambda^* + \zeta^i \cdot d^i\|^2 \\
	&= \|\lambda^{i} - \lambda^*\|^2 +2\zeta^i d^{i{\top}}(\lambda^{i} - \lambda^*) + (\zeta^i)^2\|d^i\|^2\\
	&\leq \|\lambda^{i} - \lambda^*\|^2 +2\zeta^i(D(\lambda^{i}) - D^*) + (\zeta^i)^2b^2.
	\end{align*}
	
	Then, rearranging it yields that
	\begin{align*}
	D^*-D(\lambda^{i}) \leq \frac{\|\lambda^{i} - \lambda^*\|^2}{2\zeta^i} -\frac{\|\lambda^{i+1} - \lambda^*\|^2}{2\zeta^i} + \frac{\zeta^i b^2}{2}.
	\end{align*}
	
	Summing up from $i=1$ to $k$ and noting $\zeta^i \geq \zeta^{i+1}$, it follows that
	\begin{align*}
	&\sum_{i=1}^{k}(D^*-D(\lambda^{i})) \leq -\frac{1}{2\zeta^{k+1}}\|\lambda^{k+1}-\lambda^*\| + \frac{b^2}{2} \sum_{i=1}^{k} \zeta^i \\
	&~~~+ \frac{1}{{2\zeta^1}}\|\lambda^{1} - \lambda^*\|^2 +  \frac{1}{2}\sum_{i=1}^{k}(\frac{1}{\zeta^{i+1}} - \frac{1}{\zeta^i})\|\lambda^{i+1}-\lambda^*\|^2 \\
	&\leq \frac{2}{\zeta^k}e^2 + \frac{b^2}{2} \sum_{i=1}^{k} \zeta^i.
	\end{align*}
	
	By Jenson's inequality, one can easily obtain that
	\begin{align*}
	D^*-	D(\bar{\lambda}^k) \leq \frac{2}{k\zeta^k}e^2 + \frac{b^2}{2k} \sum_{i=1}^{k} \zeta^i.
	\end{align*}
	
	The proof follows by noting that $\zeta^i = \frac{1}{be} \sqrt \frac{2}{i}$.
\end{proof}
	
By Theorem \ref{theo:model-based}, Algorithm \ref{alg} converges at a sublinear rate to an optimal policy $u(x, \lambda^*)$.
\begin{figure}[t]
	\centering
	\subfigure[The position $x_{k,1}$.]{
		\includegraphics[height=34mm]{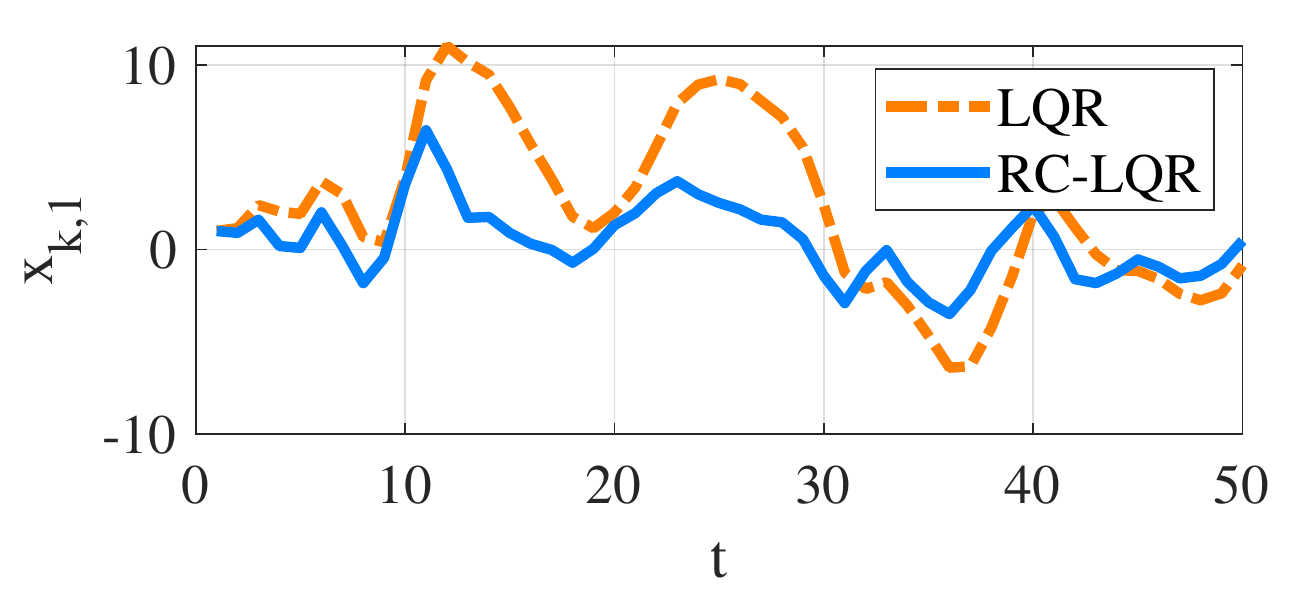}}
	
	\subfigure[The position $x_{k,3}$.]{
		\includegraphics[height=34mm]{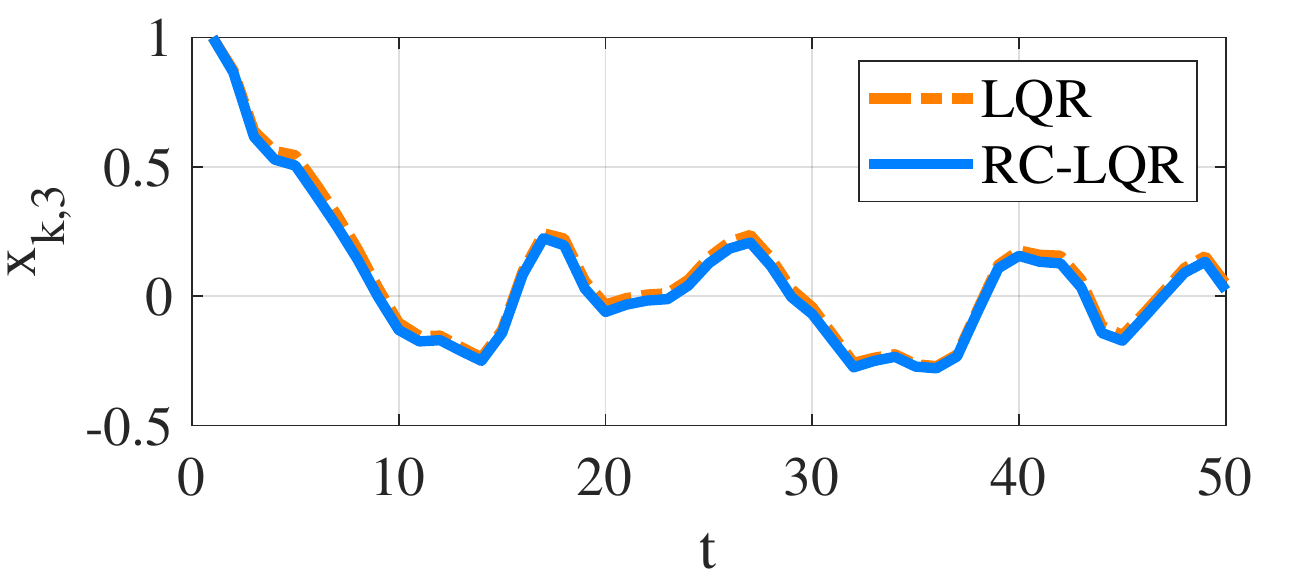}}
	\caption{The evolution of the position $x_{k,1}$ and $x_{k,3}$.}
	\label{pic:x1}
\end{figure}
\section{Simulations}\label{sec:exp}
In this section, we first demonstrate the effectiveness of our infinite-horizon risk-constrained LQR via a numerical example. Then, we validate the proposed primal-dual algorithm by examining the optimality gap and the constraint violation.
\subsection{Experimental Example}
We consider an unmanned aerial vehicle (UAV) that operates in a 2-D plane. Its discrete-time dynamical model is given by a double integrator as
\begin{equation}\label{def:model}
x_{k+1}=\begin{bmatrix}
1 & 0.5 & 0 & 0 \\
0 & 1 & 0 & 0 \\
0 & 0 & 1 & 0.5 \\
0 & 0 & 0 & 1
\end{bmatrix} x_{k}+\begin{bmatrix}
0.125 & 0 \\
0.5 & 0 \\
0 & 0.125 \\
0 & 0.5
\end{bmatrix}(u_{k}+w_{k}),
\end{equation}
where $(x_{k,1}, x_{k,3})$ is the position, $(x_{k,2}, x_{k,4})$ denotes the velocity, $u_k$ represents the acceleration and $w_k$ is the input disturbance from the wind. Suppose that the gust $w_{k,1}$ in the direction of $x_{k,1}$ is subject to a mixed Gaussian distribution of $\mathcal{N}(3,30)$ and $\mathcal{N}(8,60)$ with weights 0.2 and 0.8, respectively. In contrast, the gust $w_{k,2}$ in the orthogonal direction satisfies $w_{k,2} \sim \mathcal{N}(0,0.01)$.

We set the penalty matrix in (\ref{prob:new_rclqr}) as
$$
Q = \text{diag}(1,0.1,2,0.2)~~\text{and}~~R = \text{diag}(1,1).
$$
The risk tolerance is set to $\rho = 8$. We obtain the risk-constrained controller via Algorithm \ref{alg}. For a comparison, we compute a LQR controller, where we add an additional control input to eliminate the non-zero mean of mixed Gaussian noises $w_k$.

We demonstrate the effectiveness of our risk-constrained LQR (RC-LQR) formulation (\ref{equ:rclqr}) in Fig.~\ref{pic:x1}. It can be observed that the risk-aware controller largely compensates the risk in the state $x_{k,1}$, while its impact on the less risky state $x_{k,3}$ is consistent with that of the LQR. A more detailed discussion can be found in \cite{tsiamis2020risk}.
\begin{figure}[t]
	\centering
	\subfigure[Optimality gap $|J(u(x, \lambda^k))- J(u(x, \lambda^*))|/J(u(x, \lambda^*))$.]{
		\includegraphics[width=70mm]{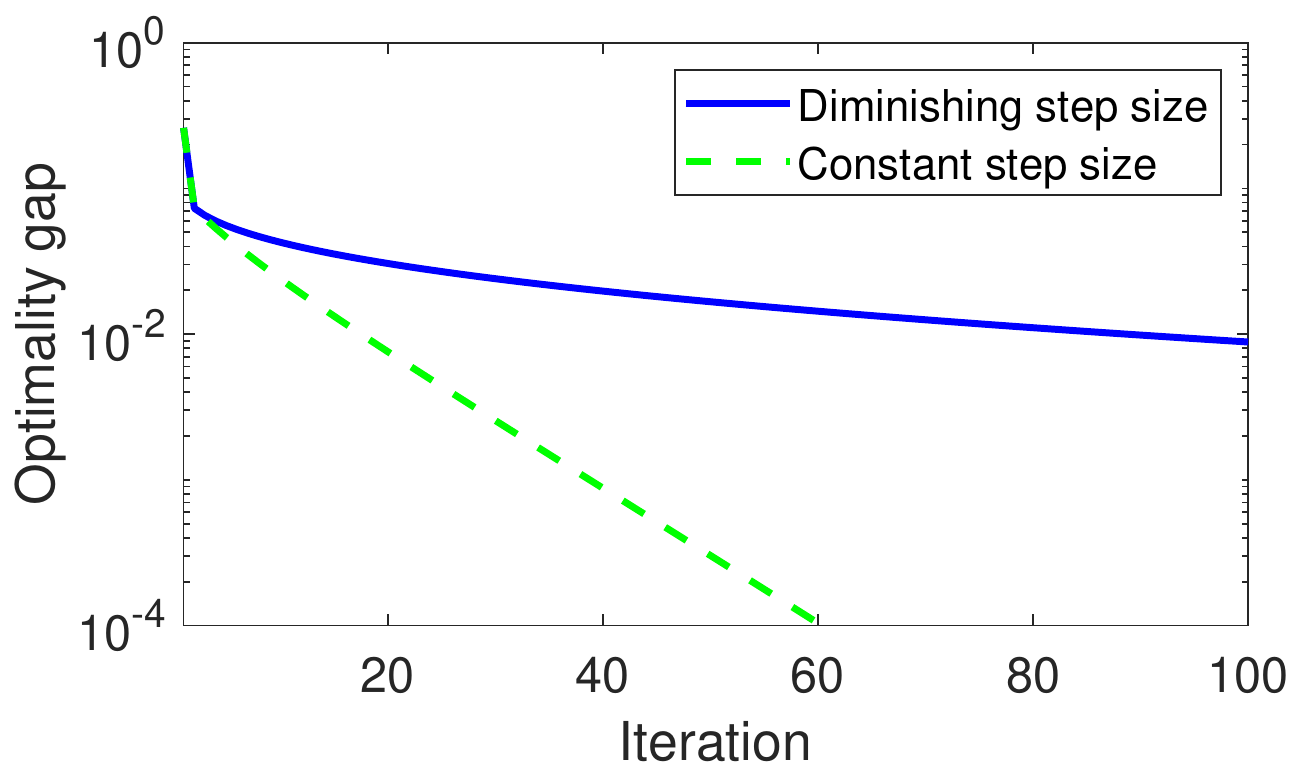}}
	
	\subfigure[Risk constraint violation $(J_c(u(x, \lambda^k))- \bar{\rho})/\bar{\rho}$.]{
		\includegraphics[width=70mm]{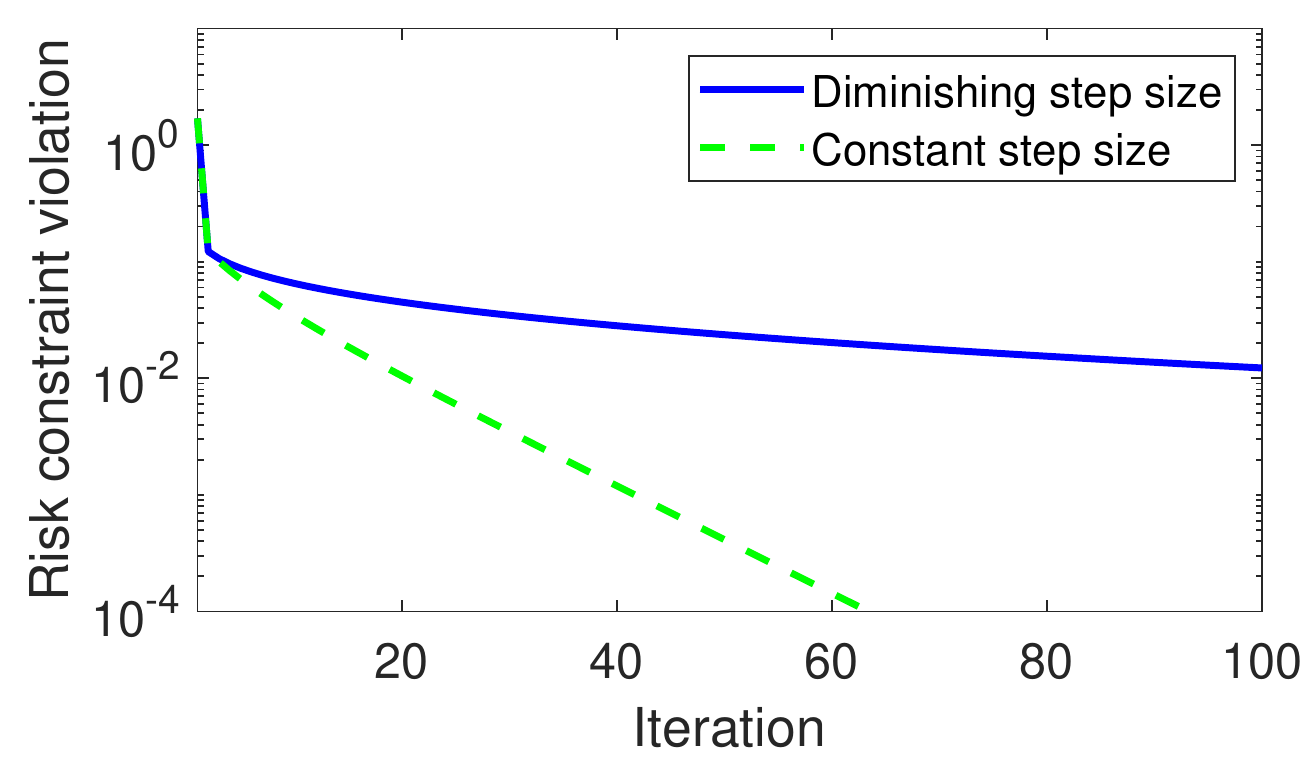}}
	\caption{Convergence of the primal-dual method.}
	\label{pic:pd_mb}
\end{figure}
\subsection{Performance of the Primal-dual Method}
We validate our primal-dual method in Algorithm \ref{alg}. We set the risk tolerance as $\bar{\rho} = 15$ in (\ref{prob:new_rclqr}), the initial multiplier as $\lambda_1 = 0$ and the diminishing step size as $\zeta^k = \frac{1}{10\sqrt{k}}$.. Since the diminishing step size rule may be overly conservative, we additionally perform Algorithm \ref{alg} with constant step size $\zeta^k = 0.1$.

Fig. \ref{pic:pd_mb} displays the optimality gap of the LQR cost and the risk constraint violation during the primal-dual optimization. The optimal cost $J^*$ is computed by $J(u(x, \lambda^*))$. Clearly, both of them converge faster under the constant step size rule. Even with the diminishing step size, the optimality gap and constraint violation reduce to less than $1\%$ within 100 iterations, exhibiting excellent performance of our model-based policy gradient primal-dual method.

\section{Conclusion}
In this paper, we have shown that an optimal policy to the infinite-horizon risk-constrained LQR problem is stationary and affine in the state. Moreover, we have proposed a primal-dual method to search an optimal policy with sublinear convergence rate.

We note that the proposed primal-dual method is model-based, i.e., the explicit dynamical model must be exactly known. Reinforcement learning, as an instance of adaptive control, has achieved tremendous success in the continuous control field. In \cite{zhao2020primal}, we have studied the model-free learning of the infinite-horizon risk-constrained LQR, which will be presented at the 3rd Annual Conference on Learning for Dynamics and Control (L4DC).

\bibliographystyle{ieeetran}
\bibliography{mybibfile}

%
%
%
%
%
%
%
%

\end{document}